\documentclass{article}
\usepackage[T1]{fontenc}                
\usepackage[utf8]{inputenc}           
\usepackage{lmodern}
\usepackage{etex}			

\usepackage[a4paper, left=20mm, right=20mm, top=10mm, bottom=15mm, headheight=10mm, headsep=6mm, footskip=10mm]{geometry}
\usepackage{fancyhdr}				
\usepackage[leftbars]{changebar}        	
\usepackage{amsmath}    					
\usepackage{amssymb}    					
\usepackage{amsfonts}  					
\usepackage{mathrsfs}   					
\usepackage{stmaryrd}  				    
\usepackage{amsthm}						
\usepackage{pifont}                      
\usepackage[all]{xy}    
\usepackage{amstext}
\usepackage{verbatim}   
\usepackage{amsbsy}
\usepackage{amsthm}
\usepackage{upgreek} 		

\usepackage{array}					   
\usepackage{supertabular}               
\usepackage{multicol}                   

\usepackage{epsfig}						
\usepackage{graphicx}					
\usepackage{rotating}					
\usepackage{pgf,tikz}					
\usetikzlibrary{arrows}

\usepackage{hyperref}  				

\usepackage{xspace}                     			
\usepackage{xcolor}								
\usepackage{calc}
\usepackage[english]{babel}						
\usepackage{indentfirst}        					
\usepackage[babel=true,kerning=true]{microtype} 	
\usepackage{numprint}							

\usepackage{datetime} 

\let\oldthebibliography=\thebibliography
\renewcommand{\thebibliography}[1]{
    \renewcommand{\bibname}{ \begin{center} \textsc{Bibliography} \end{center} }
    \oldthebibliography{\textsc{RDO}\scriptsize{4}}
    \thispagestyle{empty}
}

\newcommand{\bz}{\ensuremath{\mathbf{z}}}
\newcommand{\Pab}{\ensuremath{P^\textrm{ab}}}

\newcommand{\Hscr}{\ensuremath{\mathscr{A}}}
\newcommand{\GL}[1]{\ensuremath{\mathrm{GL}}(#1)}
\newcommand{\ZZ}{\mathbb{Z}}
\newcommand{\NN}{\mathbb{N}}
\newcommand{\UU}{\mathbb{U}}
\newcommand{\CC}{\mathbb{C}}
\newcommand{\Sfr}{\mathfrak{S}}
\newcommand{\Afr}{\mathfrak{A}}
\renewcommand{\ni}{\noindent}
\newcommand{\extensionuu}[5]{\ensuremath{\xymatrix{1 \ar[r]& #1 \ar^-{#4}[r] & #2 \ar^-{#5}[r] & #3\ar[r] & 1}}}
\newcommand{\extensionzu}[5]{\ensuremath{\xymatrix{0 \ar[r]& #1 \ar^-{#4}[r] & #2 \ar^-{#5}[r] & #3\ar[r] & 1}}}
\newcommand{\barre}[1]{\overline{#1}}
\newcommand{\wt}[1]{\widetilde{#1}}
\renewcommand{\hom}[3]{\ensuremath{\mathrm{Hom}_{#1}(#2,#3)}}
\newcommand{\Hom}[3]{\hom{#1}{#2}{#3}}
\newcommand{\ra}{\rightarrow}
\newcommand{\Cscr}{\ensuremath{\mathscr{C}}}
\newcommand{\Pscr}{\ensuremath{\mathscr{P}}}
\newcommand{\al}{\alpha}
\newcommand{\si}{\sigma}
\newcommand{\id}{\ensuremath{\mathrm{id}}}
\newcommand{\isom}[1]{\ensuremath{\overset{{\rm {\textrm{ \tiny #1}}}}{\simeq}}}
\newcommand{\diag}{\ensuremath{\mathrm{diag\,}}}
\newcommand{\email}[1]{\href{mailto:#1}{\nolinkurl{#1}}}

\newcommand{\numtoi}[1]{
	\ifthenelse{ \equal{#1}{1} }{i}{
	\ifthenelse{ \equal{#1}{2} }{ii}{
	\ifthenelse{ \equal{#1}{3} }{iii}{
	\ifthenelse{ \equal{#1}{4} }{iv}{
	\ifthenelse{ \equal{#1}{5} }{v}{
	\ifthenelse{ \equal{#1}{6} }{vi}{
	\ifthenelse{ \equal{#1}{7} }{vii}{
	\ifthenelse{ \equal{#1}{8} }{viii}{
	\ifthenelse{ \equal{#1}{9} }{ix}{
	\ifthenelse{ \equal{#1}{10} }{x}{
	\ifthenelse{ \equal{#1}{11} }{xi}{
	\ifthenelse{ \equal{#1}{12} }{xii}{
	\ifthenelse{ \equal{#1}{13} }{xii}{
	\ifthenelse{ \equal{#1}{14} }{xiv}{
	\ifthenelse{ \equal{#1}{15} }{xv}{
	\ifthenelse{ \equal{#1}{16} }{xvi}{
	\ifthenelse{ \equal{#1}{17} }{xvii}{
		ERREUR,~MODIFIER~LA~MACRO~numtoi
	}	}	}	}	}	}	}	}	} } } } } } } } }
}	

\newcounter{cretraiti}
\newenvironment{ri}[1]{\begin{list}{($\numtoi{\thecretraiti}$)}{
	\usecounter{cretraiti}
	\topsep=0.5ex
 	\itemsep=0.3ex
 	\labelsep=0.3em
 	\parsep=0ex
 	\listparindent=1em
 	\settowidth{\labelwidth}{(#1)}
 	\leftmargin=\labelwidth
}}{\end{list}
}

\newcommand{\spacebeforeenv}{\vspace{1ex}}
\newcommand{\spaceafterenv}{\vspace{1ex}}

\newcommand{\tiret}{\rule[0.6ex]{1.3ex}{0.26ex}}

\newcommand{\myqedenv}{}

\newcommand{\envfont}{\sf\bfseries}

\newcounter{ctheo}
\renewcommand{\thectheo}{\arabic{ctheo}}

\newcommand{\metaenvironnementhm}[2]  
{
  \refstepcounter{ctheo}
  \ifthenelse{ \equal{#2}{toto} }{
    \spacebeforeenv \begin{traitsurlecote} \ni {\envfont #1 \thectheo}
  }{
    \spacebeforeenv \begin{traitsurlecote} \ni {\envfont #1 \thectheo{} \tiret\, #2.}
  }
}

\newenvironment{traitsurlecote}
{\cbstart
\setcounter{changebargrey}{0}      
}
{\cbend
}

\newenvironment{theo}[1][toto]
{
  \metaenvironnementhm{Theorem}{#1}
}
{\end{traitsurlecote}\spaceafterenv}

\newenvironment{prop}[1][toto]
{
  \metaenvironnementhm{Proposition}{#1}
}
{\end{traitsurlecote}\spaceafterenv}

\newenvironment{lem}[1][toto]
{ \metaenvironnementhm{Lemma}{#1} }
{\end{traitsurlecote}\spaceafterenv}

\newenvironment{cor}[1][toto]
{ \metaenvironnementhm{Corollary}{#1} }
{\end{traitsurlecote}\spaceafterenv}

\newenvironment{prop-dfn}[1][toto]
{ \metaenvironnementhm{Proposition-Definition}{#1} }
{\end{traitsurlecote}\spaceafterenv}

\newenvironment{demo}
{\ni  {\envfont Proof.~}}
{\spaceafterenv}

\newcommand{\metaenvironnement}[2]  
{
  \refstepcounter{ctheo}
  \ifthenelse{ \equal{#2}{toto} }{
    \spacebeforeenv \ni {\envfont #1 \thectheo}
  }{
    \spacebeforeenv \ni {\envfont  #1 \thectheo{} \tiret\, #2.}
  }
}

\newcommand{\metaenvironnementsimple}[2]  
{
  \ifthenelse{ \equal{#2}{toto} }{
    \spacebeforeenv \ni {\envfont #1}      
  }{
    \spacebeforeenv \ni {\envfont #1 \tiret\, #2.}
  }
}

\newenvironment{ex}[1][toto]
{ \metaenvironnement{Example}{#1} }
{\myqedenv\spaceafterenv}

\newenvironment{rem}[1][toto]
{ \metaenvironnement{Remark}{#1} }
{\myqedenv\spaceafterenv}

\begin{document}

\begin{center}
{\LARGE\bf Torsion subgroups of quasi-abelianized braid groups}
\end{center}

\hspace*{1cm}
\begin{minipage}{6cm}
Vincent Beck
\newline MAPMO (UMR CNRS 7349),
\newline Université d'Orl\'eans, 
\newline F-45067, Orl\'eans, France 
\newline FDP - FR CNRS 2964
\newline
\email{vincent.beck@univ-orleans.fr}
\end{minipage}
\hspace*{2cm}
\begin{minipage}{6cm}
Ivan Marin
\newline LAMFA (UMR CNRS 7352),
\newline Université de Picardie Jules Verne, 
\newline 80039 Amiens Cedex 1 France 
\newline
\email{ivan.marin@u-picardie.fr}
\end{minipage}

\begin{abstract} This article extends the works of Gonçalves, Guaschi, Ocampo~\cite{ggo} 
and Marin~\cite{marin2} on finite subgroups of the quotients of generalized braid groups by the derived subgroup 
of their pure braid group. We get explicit criteria for subgroups of the (complex) reflection group
to lift to subgroups of this quotient. In the specific case of the classical braid group, this enables us to
describe all its finite subgroups : we show that every odd-order finite group can be embedded
in it, when the number of strands goes to infinity. We also determine a complete list of the irreducible
reflection groups for which this quotient is a Bieberbach group.
\end{abstract}

\section{Introduction}

Let $V$ be a finite dimensional complex vector space and $W \subset \GL{V}$ be an irreducible complex reflection group.  We define $\Hscr$ to be the set of reflection hyperplanes for $W$ and $V^\textrm{reg}=V \setminus \cup_{H \in \Hscr} H$. 
Let $x_0 \in V^\textrm{reg}$. We consider $P=\pi_1(V^\textrm{reg},x_0)$ the pure braid group of $W$ and 
$B=\pi_1(V^\textrm{reg}/W,\overline{x}_0)$ its braid group where $\overline{x}_0$ is the image of $x_0$ in $V^\textrm{reg}/W$. 
We have the following defining exact sequence (see~\cite{bmr} for more on this subject)
$$\extensionuu{P}{B}{W}{}{}$$
\ni giving rise to the abelian extension of $W$ 
$$\extensionzu{\Pab}{B/[P,P]}{W}{i}{p}$$
\ni This last extension whose order has been studied in~\cite{beck-pbw} will play a key role in the sequel. 
Let us say that $B/[P,P]$ is the {\it relative abelianization of the braid group of $W$}. Let us also remind that as a $W$-module 
$\Pab$ is the permutation module $\ZZ \Hscr$ (see~\cite[Section 6.1]{orlik-terao}).

We get a very concrete criterion for a subgroup of $W$ to embed in $B/[P,P]$ in terms of the stabilizers of the hyperplanes 
of $W$ (Theorem~\ref{theo-finite-order}); for the case of the symmetric group $\Sfr_n$ (and more generally for a wide class
of complex reflection groups including the real ones) and the standard braid group $B_n$, we are able to 
describe all finite subgroups of $B/[P,P]$ (see Corollary~\ref{cor-real-reflection-group}
and~Proposition~\ref{prop-subgroup-binfty}). In particular we show that every odd-order subgroup of $\Sfr_n$ 
can be embedded in $B_n/[P_n,P_n]$ and that every odd-order group can be embedded in $B_\infty/[P_\infty,P_\infty]$ where $B_\infty$
(resp. $P_\infty$) denotes the infinite (pure) braid group. We are also able to give a complete list of the irreducible reflection groups such that $B/[P,P]$ is a Bieberbach group (Corollary~\ref{cor-bieberbach}).

Let us start with the notations we will use in this article. 
The group $W$ acts on $\Hscr$: if $H$ is the hyperplane of the reflection $s \in W$, 
then $wH$ is the hyperplane of the reflection $wsw^{-1}$. Moreover $wH=H$ if and only if $wsw^{-1}=s$. 
For $H \in \Hscr$, we denote by $N_H=\{w \in W, wH=H\}= \{w \in W, sw=ws\}$ the stabilizer of $H$ which is also the 
centralizer of any reflection of $W$ with hyperplane $H$.  

For a subgroup $G$ of $W$ and $H \in \Hscr$, the stabilizer of $H$ under the action of $G$ is $N_H\cap G$. 
We denote it by $N_{H,G}$. 

Since $W$ is generated by reflections, we then deduce  $w H=H$ for every $H \in \Hscr$ if and only $w \in Z(W)$ where 
$Z(W)$ is the centre of $W$. In particular, the group $\barre{W}=W/Z(W)$ acts faithfully on $\Hscr$. 

All along this article, we will consider subgroups $G$ of $W$ and subgroups $\barre{G}$ of $\barre{W}$. 
For a subgroup $G$ of $W$, we denote by $\barre{G}$ its direct image in $\barre{W}$. 

Finally since $W$ is finite, we may assume that $\langle \cdot, \cdot \rangle$ is an hermitian product invariant under $W$. 
For $H \in \Hscr$, we set $C_H=\{w \in W,\ \forall\,x \in H^\perp, \ w(x)=x\}$ the parabolic subgroup of $W$ associated to $H^\perp$.

Let us end this introduction with a cohomological lemma that we will need later. 

\begin{lem}\label{lem-ext} 
Let $G$ be a finite group, $X$ a set on which $G$ acts and $\ZZ X$ the corresponding permutation module. 
Then $H^1(G,\ZZ X)$ is trivial. 
\end{lem}

\begin{demo} Decomposing $X$ into $G$-orbits $X=\sqcup_{A \in X/G} A$, we get $H^1(G,\ZZ X) = \bigoplus_{A \in X/G} H^1(G,\ZZ A)$. 
Choosing for each $A \in X/G$ a representative $a$ and using Shapiro's isomorphism (\cite[Proposition III.6.2]{brown}) 
we obtain that $H^1(G,\ZZ X) = \bigoplus_{A \in X/G} H^1(G_a,\ZZ)$ where $G_a$ is the stabilizer of $a$. 
But $G_a$ acts trivially on $\ZZ$ and then $H^1(G_a,\ZZ)=\Hom{\textrm{gr.}}{G_a}{\ZZ}=\{0\}$ since $G_a \subset G$ is finite. 
\end{demo}

\section{Finite subgroups of $B/[P,P]$}

\subsection{A criterion for finite subgroup}

In this subsection, we give a criterion for a subgroup of $W$ to lift in $B/[P,P]$. The criterion~(Theorem~\ref{theo-finite-order})
has many consequences such that the description of the finite subgroups of $B_n/[P_n,P_n]$ and $B_\infty/[B_\infty,B_\infty]$
which are explored in the following subsections.  

\begin{theo}\label{theo-finite-order} Let $G$ be a subgroup of $W$ then the extension 
$$\extensionuu{\Pab}{p^{-1}(G)}{G}{i}{p}$$
\ni splits if and only if for every $H \in \Hscr$ we have $N_H \cap G \subset C_H$.
When these conditions are fulfilled, $G$ identifies to a finite subgroup of $B/[P,P]$ and $\barre{G}$ is isomorphic to $G$.
Moreover $p^{-1}(G)$ is a crystallographic group. 

Let $\wt{G}$ be a torsion subgroup of $B/[P,P]$ then $\wt{G}$ is finite and is isomorphic to its image in $W$ and the extension 
$$\extensionuu{\Pab}{p^{-1}(p(\wt{G}))}{p(\wt{G})}{i}{p}$$
\ni splits. Moreover every subgroup of $B/[P,P]$ isomorphic through $p$ to $p(\wt{G})$ is conjugated to $\wt{G}$ 
by an element of $\Pab$ and the intersection of the normalizer of $\wt{G}$ in $B/[P,P]$ with $\Pab$ is isomorphic to $\ZZ \Hscr/G$. 
\end{theo}

\begin{demo} Let us consider the extension $\extensionuu{\Pab}{p^{-1}(G)}{G}{i}{p}$ as an element in $H^2(G,\Pab)$. 
This is the image of the extension $\extensionuu{\Pab}{B/[P,P]}{W}{i}{p}$ through the restriction map 
$\textrm{Res}\colon H^{2}(W,\Pab) \ra H^2(G, \Pab)$.  

Following the proof of Proposition~5 of~\cite{beck-pbw}, we get the following commutative diagram
$$\displaystyle\xymatrix{H^2(W,\Pab) \ar[r] \ar^{\textrm{Res}}[d] & \left[\bigoplus_{H \in \Hscr} H^2(N_H,\ZZ)\right]^W \ar[r]
\ar^{\textrm{Res}}[d]& \left[\bigoplus_{H \in \Hscr} \Hom{\textrm{gr}}{N_H}{\CC^\times} \right]^W \ar^{\textrm{Res}}[d]\\
	H^2(G,\Pab) \ar[r] & \left[\bigoplus_{H \in \Hscr} H^2(N_H\cap G,\ZZ)\right]^G  \ar[r]& 
	\left[\bigoplus_{H \in \Hscr} \Hom{\textrm{gr}}{N_H \cap G}{\CC^\times} \right]^G }$$
\ni Hence the extension~$\extensionuu{\Pab}{p^{-1}(G)}{G}{i}{p}$ can be identified with the family of group homomorphisms 
$r_H: N_H \cap G \ra \CC^\times$ given by the restriction to the line $H^\perp$ whose kernel is by definition $C_H$. 

In this case, $G \cap ZW=\{1\}$ since for every non trivial element $w$ of the center of $W$ and every hyperplane $H \in \Hscr$,
we have $r_H(w) \neq 1$. So $\barre{G}$ and $G$ are isomorphic. 
Finally, $p^{-1}(G)$ is a crystallographic group since $G$ acts faithfully on $\Hscr$ since $G \cap ZW=\{1\}$.

Let us now consider $\wt{G}$ a  torsion subgroup of $B/[P,P]$. For $x \in \wt{G}$, let us first show that the order of $x$
is the same that the order of $p(x)\in W$. Let $m$ be the order of $x$ and assume that $p(x)^n=1$ for some $n$ dividing $m$
and $n \neq m$. Then $x^n \in P$ but $x^n \notin [P,P]$. But $\Pab$ is a torsion-free group, so it is impossible. 
The map $p:\wt{G} \ra W$ is then injective and $\wt{G}$ is finite and $p: \wt{G} \ra p(\wt{G})$ is an isomorphism.

The $\Pab$-conjugacy classes of subgroups of $B/[P,P]$ isomorphic to $p(\wt{G})$ through $p$ are in bijection with 
$H^{1}(p(\wt{G}),\Pab)$ (see~\cite[Proposition 2.3 Chapter IV]{brown}). But $H^1(p(\wt{G}),\Pab)$ is trivial (Lemma~\ref{lem-ext}) 
and  we get the result.  

Finally, let $x \in \Pab$ such that $x\wt{G}x^{-1}=\wt{G}$. Then for every $g \in \wt{G}$, we have $xgx^{-1}g^{-1} \in 
\wt{G} \cap \Pab$. But $\wt{G}$ is finite and $\Pab$ is torsion-free, hence $xgx^{-1}g^{-1}=1$. 
We then deduce $gxg^{-1}=x$ for all $g \in \wt{G}$ and so 
$$\displaystyle x= \sum_{\Cscr \in \Hscr/p(\wt{G})} \al_\Cscr \left(\sum_{H \in \Cscr} c_H\right)\,.$$
\end{demo}

Since $w \in W$ has a finite order lifting in $B/[P,P]$ if and only if 
the extension
$$\extensionuu{\Pab}{p^{-1}(\langle w \rangle)}{\langle w \rangle}{i}{p}$$
\ni is split, we get the following corollary.

\begin{cor}\label{cor-lifting-element} An element $w \in W$ has a lifting in $B/[P,P]$ of finite order if and only if 
$\langle w \rangle \cap N_H \subset C_H$ for every $H \in \Hscr$. 
\end{cor}

We can then easily generalize Theorem 2.5 of~\cite{marin2}.

\begin{cor}\label{cor-lifting} Let $w \in ZW \setminus \{1\}$; the order of every lifting of $w$ in $B/[P,P]$ is infinite.

Let $w \in W$ such that $w$ has two eigenvalues one of those is $1$ then the order of every lifting of $w$ in $B/[P,P]$ is infinite.
In particular, if $w \in W$ has order $2$ or is a reflection then the order of every lifting 
of $w$ in $B/[P,P]$ is infinite. 
\end{cor}

\begin{demo} Let $w \in ZW \setminus \{1\}$ then $w \in N_H$ for every $H \in \Hscr$ and $w \notin C_H$. 
So~Theorem~\ref{theo-finite-order} gives the result. This result follows also directly from the fact that,
when $W$ is irreducible $p^{-1}(ZW)$ is a free abelian group of rank $|\Hscr|$ generated by $\Pab$ and 
the class $\underline{\bz}$ of the path $t \mapsto \exp(2i\pi t/|Z(W)|)x_0$ (see~\cite{marin2}). 

Let us now consider $w \in W$ with $w$ having two eigenvalues $1$ and $\lambda \neq 1$. 
Steinberg's theorem implies that there exists a reflection in $W$ such that its hyperplane $H$ contains the eigenspace of $w$
associated to $1$. Then $H^\perp$ is contained in $\ker(w-\lambda \id)=\ker(w -\id)^\perp$. Thus the restriction of 
$w$ on $H^\perp$ is $\lambda \id \neq \id$. 

If $w \in W$ has order $2$ then either $w = -\id \in ZW \setminus \{1\}$ or $w$ has two eigenvalues $1$ and $-1$. 
If $w \in W$ is a reflection, its eigenvalues are $1$ and $\lambda$ for some $\lambda \neq 1$. 
\end{demo}

\begin{rem} In fact, Theorem~2.5 of~\cite{marin2} and its proof can easily be extended to the case of $w \in W$ 
such that $w$ has two eigenvalues one of those is $1$. So the previous corollary was already known but we give here 
a different proof. 
\end{rem}

\begin{lem}\label{cor-power-element} If $w \in W$ is such that a power of $w$ has only infinite lifting in $B/[P,P]$ 
then every lifting of $w$ in $B/[P,P]$ has infinite order. 

This is in particular the case when $w$ is of even order or if a power of $w$ is a non trivial element of $ZW$. 
\end{lem}

\begin{demo} If $x$ is a lifting of $w$ then $x^n$ is a lifting of $w^n$.
\end{demo}

We may interpret Theorem~\ref{theo-finite-order} as a kind of local property in the following sense: a subgroup $G$
of $W$ identifies to a subgroup of $B/[P,P]$ if and only if every element of $G$ has a finite order lifting in $B/[P,P]$.
More precisely, we get the following corollary. 

\begin{cor} Let $G$ be a subgroup of $W$. The extension $\extensionuu{\Pab}{p^{-1}(G)}{G}{i}{p}$ is split
if and only if every element of $G$ has a finite order lifting in $B/[P,P]$. 
\end{cor}

\begin{demo} Both conditions are equivalent to the following: for every $H \in \Hscr$ and every 
$w \in G \cap N_H$, $w \in C_H$. 
\end{demo}

\subsection{Examples and Applications}

We start this subsection by describing the elements of the groups in the infinite series that can be lifted 
in their restricted braid group. We then study the groups $W$ such that $B/[P,P]$ is a Bieberbach group and 
finally we give a list of groups $W$ (containing the Coxeter groups) such that every odd-order element of $W$ has a finite 
order lifting in $B/[P,P]$.

For an integer $n$, we denote by $\UU_n$ the set on $n^\textrm{th}$ root of unity in $\CC$. 
Let us consider $(e_1,\ldots, e_r)$ the canonical basis of $\CC^r$. 
For $w \in G(de,e,r)$, we denote by $\si_w \in \Sfr$ the image of $w$ through the natural homomorphism 
$G(de,e,r) \ra \Sfr_r$. Thus, there exists a family $(a_i)_{1 \leq i \leq r} \in {\UU_{de}}^r$ such that $w(e_i)=a_i e_{\si_{w}(i)}$
for all $i \in \{1,\ldots, r\}$.  Let us write $\si_w=\si_1\cdots \si_s$ as a product of disjoint cycles. 
For such a cycle $\si=(i_1,\ldots, i_\al)$,
we denote by $p_{w,\si}=a_{i_1}\cdots a_{i_\al} \in \UU_{de}$.

\begin{cor}[Finite order lifting for the infinite series]\label{cor-infinite-series} Let $w \in G(de,e,r)$. 

If $w$ has a finite order lifting in $B/[P,P]$ then either $\si_w \neq \id$ and for every cycle $\si$ appearing in $\si_w$, 
we have $p_{w,\si}=1$, or $\si_w =\id$ and for every $i \neq j \in \{1,\ldots,r\}$ the order of $a_i{a_j}^{-1}$ 
is a multiple of the order of $a_i$.

If $w$ is of odd order and for every cycle $\si$ appearing in $\si_w$, we have $p_{w,\si}=1$ then $w$
has a finite order lifting in $B/[P,P]$. 

When $d \geq 2$, then $w$ has a finite order lifting in $B/[P,P]$ if and only if 
$w$ is of odd order and for every cycle $\si$ appearing in $\si_w$, we have $p_{w,\si}=1$.

When $d=1$, $w$ has a finite order lifting in $B/[P,P]$ if and only if $w$ is of odd order 
and for $w$ for every cycle $\si$ of $\si_w$, we have $p_{w,\si}=1$ or 
$\si_w =\id$ and for every $i \neq j \in \{1,\ldots,r\}$ the order of $a_i{a_j}^{-1}$ is a multiple of the order of $a_i$.
\end{cor}

\begin{demo} The result is clear for $r=1$. Let us assume that $r \geq 2$. 

Let us consider $w \in W$ with a finite order lifting in $B/[P,P]$ and $\si$ a cycle of $\si_w$ of length greater than $1$.
Set $i \neq j$ in the support of $\si$ and $\zeta \in \UU_{de}$ and consider the hyperplane 
$$H_{i,j,\zeta}=\{(z_1,\ldots,z_r) \in \CC^r, \  z_i=\zeta z_j\}.$$
\ni For $\ell \in \ZZ$, if $w^\ell \in N_{H_{i,j,\zeta}}$ then $\{{\si_w}^\ell(i),{\si_w}^\ell(j)\}=\{i,j\}$. 
But $\si_w$ is of odd order (Corollary~\ref{cor-lifting}), hence ${\si_w}^\ell(i)=i$ and ${\si_w}^\ell(j)=j$. 
We then deduce that $\ell$ is a multiple of the length $\al$ of $\si$: we can write $\ell = \alpha m$.
Then the restriction of $w^\ell$ to the space
$\textrm{Span}(e_{i_1},\ldots, e_{i_\alpha})$ (where $\si =(i_1,\ldots, i_\alpha)$) is ${p_{w,\si}}^m \id$. 
But from Corollary~\ref{cor-lifting-element}, we get that $w^\ell \in C_{H_{i,j,\zeta}}$, hence ${p_{w,\si}}^m=1$. 
Since we can choose $m=1$, we get the result.    

If $i$ is a fixed point of $\si_w$. We have $w(e_i)=a_ie_i$. We want to show that $a_i=1$ when 
there exists a non trivial orbit $\si$ of $\si_w$, we consider $j$ in the support of $\si$. 
If $w^\ell$ in $N_{H_{i,j,1}}$ then $\textrm{Span}(e_i,e_j)$ is stable by $w^\ell$. Hence $\ell$
is a multiple of the length of $\si$. Hence $w^\ell(e_j)=e_j$ by the preceeding point. But $w^\ell(e_j- e_i)=e_j- e_i$
since $w^\ell \in C_{H_{i,j,1}}$ (Corollary~\ref{cor-lifting-element}). Hence $w^\ell(e_i)=e_i$. 

Let us now consider the case where $\si_w=\id$ 
then we can write $w(e_i)=a_ie_i$ for all $i$. For $i\neq j$, we have $w^\ell \in N_{H_{i,j,1}}$ 
if and only if ${a_i}^\ell={a_j}^\ell$. 
If $w^\ell \in N_{H_{i,j,1}}$ then $w^\ell(e_i-e_j)=(e_i -e_j)$ (since $w^\ell \in C_{H_{i,j,\ell}}$). 
Hence ${a_i}^\ell={a_j}^\ell$ implies ${a_i}^\ell={a_j}^\ell=1$. So for every $i \neq j$ the order of $a_i{a_j}^{-1}$
is a multiple of the order of $a_i$ and $a_j$. 

Conversely, assume that $p_{w,\si}=1$ for every cycle $\si$ of $\si_w$ and $w$ is of odd order
and let us apply the criterion of Theorem~\ref{theo-finite-order} to show that $w$ has a finite order lifting in $B/[P,P]$. 
Let $1 \leq i \leq r$ and $H_i=\{(z_1,\ldots, z_r) \in \CC^r, \ z_i= 0\}$. 
If $w^\ell \in N_{H_i}$ then $w^\ell(e_i) \in \CC^\times e_i$. 
Hence $\ell$ is a multiple of the length of the cycle $y=(i,i_1,\ldots, i_\al)$ of $\si_w$ whose support contains $i$.
So that the restriction of $w^\ell$ to $\textrm{Span}(e_i,e_{i_1},\ldots, e_{i_\alpha})$ is the identity since $p_{w,\si}=1$.
So $w^\ell \in C_{H_i}$. 

For $i \neq j \in \{1,\ldots, r\}$ and $\zeta \in \UU_{de}$, if $w^\ell \in N_{H_{i,j,\zeta}}$ 
then $\{{\si_w}^\ell(i),{\si_w}^\ell(j)\} = \{i,j\}$. 
If $i$ and $j$ belongs to different orbits $\si$ and $\si'$ under $\si_w$ then necessarily ${\si_w}^\ell(i)=i$ 
and ${\si_w}^\ell(j)=j$.  Hence  $\ell$ is a multiple of the lengths of $\si$ and $\si'$. But $p_{w,\si}=p_{w,\si'}=1$ 
so that $w^\ell(e_i)=e_i$ and $w^\ell(e_j)=e_j$ and $w^\ell \in C_{H_{i,j,\zeta}}$. If $i$ and $j$ are in the 
same orbit $\si$ under $\si_w$ then we also have
${\si_w}^\ell(i)=i$ and ${\si_w}^\ell(j)=j$ because if ${\si_w}^\ell(i)=j$ and ${\si_w}^\ell(j)=i$ then $\si_w$ would be 
of even order and $w$ too, which is not. Hence $\ell$ is a multiple of the length of $\si$ and since $p_{w,\si}=1$, we 
get that $w^\ell \in C_{i,j,\zeta}$.

Let us now consider the case $d\geq 2$. From the first two parts, the only thing to show is that if 
$w$ verifies $\si_w = \id$ and $w$ has finite order lifting in $B/[P,P]$ then $w=\id$.  
Let $i \in \{1,\ldots,r\}$. We have $w(e_i)=a_i e_i$. Hence $w \in N_{H_i}$ 
where $H_i=\{(z_1,\ldots, z_r) \in \CC^r, \ z_i= 0\}$. We then deduce that $w \in C_{H_i}$ that is $a_i=1$. 

Let us consider the case $d=1$. From the first two parts, the only thing to show is that if $\si_w=\id$
and for $i\neq j \in \{1,\ldots,r\}$, the order of $a_i{a_j}^{-1}$ is a multiple of the order of $a_i$ then 
$w$ has a finite order lifting in $B/[P,P]$. 
For $i \neq j \in \{1,\ldots, r\}$ and $\zeta \in \UU_{de}$, if $w^\ell \in N_{H_{i,j,\zeta}}$ 
then $w^\ell$ stabilizes three lines in the plane $\textrm{Span}(e_i,e_j)$; namely $\CC e_i,\CC e_j$ and 
the line ${H_{i,j,\zeta}}^\perp$. So ${a_i}^\ell={a_j}^\ell$ and then ${a_i{a_j}^{-1}}^\ell=1$. 
Hence ${a_i}^\ell= {a_j}^\ell= 1$ and then $w^\ell \in C_{H_{i,j,\zeta}}$. Corollary~\ref{cor-lifting-element}
gives then the result. 
\end{demo}


\begin{ex}\label{ex-diag-matrix} 
Let $j = -1/2+ i \sqrt{3}/2$. The preceding corollary  ensures us that $w=\diag(j,j^2,1\ldots, 1) \in G(3,3,r)$ 
verifies $\si_w = \id$ and has a finite order lifting in $B/[P,P]$. In particular, for $r=2$, $\diag(j,j^2)$
gives an example of an element of $G(3,3,2)$ of order $3$ that has a finite order lifting in $B/[P,P]$ but such that no element of 
order $3$ of $\Sfr_2$ (since there are no such element) lifts in $B/[P,P]$. 

More generally, if $de = 2^m q$ with $q$ odd, let $\zeta$ of $q^\textrm{th}$ root of unity. Then 
$w=\diag(\zeta,\zeta^{-1},1,\ldots, 1) \in G(de,de,r)$ has a finite order lifting in $B/[P,P]$. 
\end{ex}

\begin{rem} Let $G$ be a subgroup of $G(de,e,r)$ where $d \geq 2$ which lifts in $B/[P,P]$. 
Then Corollary~\ref{cor-infinite-series} shows that $G \cap D=\{1\}$ where $D$ is the subgroup of diagonal matrices 
of $G(de,e,r)$. In particular $G$ is isomorphic to its image in $\Sfr_r$. 
Moreover the criterion of Corollary~\ref{cor-infinite-series}
shows also that every odd-order subgroup of $\Sfr_r$ lifts into $B/[P,P]$. 
So that the isomorphism classes of finite subgroups of $B/[P,P]$ are the same that the isomorphism classes of odd-order 
subgroups of $\Sfr_r$. The preceding example shows that this is not the case when $d=1$. 
\end{rem}

In~\cite{marin2}, Marin shows that $B/[P,P]$ is a Bieberbach group for $G_4$ and $G_6$. 
We are able to give the list of the irreducible group such  that $B/[P,P]$ is a Bieberbach group (Corollary~\ref{cor-bieberbach}).
But we start with two small groups that we can study by hands.

\begin{cor} For $W=G_5$, the group $B/[P,P]$ is a Bieberbach group of holonomy group $\Afr_4$ and dimension~$8$.

For $W=G_7$, the group $B/[P,P]$ is a Bierberbach group of holonomy group $\Afr_4$ and dimension~$14$. 
\end{cor}

\begin{demo} In $G_5$, there are elements of order $1,2,3,4,6,12$. 
Corollary~\ref{cor-power-element} ensures us that we have only to study the case of the elements of order $3$. 
There are $26$ of them shared in 8 conjugacy classes. There are two elements in the center of $G_5$. 
There are $16$ reflections. Both cases are treated by previous results. It remains to consider two conjugacy classes which 
are products of an element of order $3$ of $ZG_5$ with a reflection and are not reflections: they have two eigenvalues which 
are $j$ and $j^2$ (the third non trivial roots of unity). Such an element is contained in the stabilizer of a hyperplane
but not in the corresponding $C_H$. So Theorem~\ref{theo-finite-order} allows us to conclude. 

The situation is exactly the same for $G_7$ since the elements of $G_7$ have order $1,2,3,4,6,12$ and 
every element of order $3$ of $G_7$ is contained in the $G_5$ subgroup of $G_7$.
\end{demo}

\begin{cor}[Bieberbach groups]\label{cor-bieberbach} For $W=G(de,e,r)$, then $B/[P,P]$ is a Bieberbach group if and only if
$r=1$ or $r=2$ and $d \geq 2$ or $r=2$, $d=1$ and $e=2^m$ for some $m$. 

Among the exceptionnal groups, the group $B/[P,P]$ is a Bieberbach group for 
$$W=G_4, G_5, G_6, G_7, G_{10}, G_{11}, G_{14}, G_{15}, G_{18}, G_{19}, G_{25}, G_{26}.$$ 
\end{cor}

\begin{demo} Assume that for $G(de,e,r)$, $B/[P,P]$ is a Bieberbach group. 
Corollary~\ref{cor-infinite-series} ensures us that $(1,2,3)\in \Sfr_r$ has a finite order lifting in $B/[P,P]$. 
Hence $r \leq 2$. If $r=2$ and $d=1$, then~Example~\ref{ex-diag-matrix} shows that $e$ has to be a power of $2$. 

Conversely, if $r=1$, then $B=B/[P,P]=\ZZ$. If $r=2$ and $d \geq 2$ then a diagonal matrix $D \in G(de,e,2)$ 
has a finite order lifting only if $D=\id$ (since for every $i\in \{1,2\}$ we have $a_i=p_{D,\{i\}}=1$). 
And a non diagonal matrix is of even order so it does not have a finite order lifting in $B/[P,P]$. 

For $r=2,d=1,e=2^m$, let $D=\diag(\zeta,\zeta^{-1})$ be a diagonal matrix of $G(e,e,2)$ with $\zeta \in \UU_e$.
If $\zeta\neq 1$ then the order of $\zeta^2$ is half the order of $\zeta$ (since $e$ is a power of $2$) and 
by Corollary~\ref{cor-infinite-series}, $D$ does not have a finite order lifting in $B/[P,P]$.  
A non diagonal matrix is of even order so it does not have a finite order lifting in $B/[P,P]$.

For the exceptionnal groups, we use the package \cite{chevie} of~\cite{gap} to compute stabilizers of hyperplanes and 
parabolic subgroups associated to the orthogonal line.
\end{demo}

We now give a corollary generalizing Corollary 2.7 of~\cite{marin2} and Theorem~16 of~\cite{ggo}. 

\begin{cor}[Real reflection groups]\label{cor-real-reflection-group} 
Let $W$ be a real reflection group. Then every odd-order element has a finite order lifting
in $B/[P,P]$. More generally, if $W$ is a reflection group for which the roots of unity inside its field of definition all 
have for order of power of $2$ (The irreducible groups verifying this property 
are the following $G_8,G_9,G_{12},G_{13},G_{22}, G_{23}, G_{24}, G_{28},G_{29},G_{30},G_{31},G_{35},G_{36},G_{37}$ 
for the exceptional ones, and $G(de,e,r)$ with $d$ and $e$ powers of $2$ and $r \geq 2$, $\Sfr_r$, $G(d,1,1)$
with $d$ power of $2$ and $G(e,e,2)$ with $e \geq 3$ for the infinite series), 
then every odd-order element has a finite order lifting in $B/[P,P]$.

For those groups, let $G \subset W$ then $p^{-1}(G)$ the inverse image of $G$ in $B/[P,P]$ is a Bieberbach group 
if and only if $G$ is a $2$-subgroup of $W$. 

Again, for those groups, let $G \subset W$ with $|G|$ odd then $B/[P,P]$ has a finite group isomorphic to $G$
and $p^{-1}(G)$ is a crystallographic group.
\end{cor}

\begin{demo} Let $w$ be an odd-order element of $W$ and $H \in \Hscr$. If $w \in N_H$ then $w$ stabilizes $H^\perp$
and the restriction of $w$ to $H^\perp$ is still of odd order. In particular the eigenvalue of $w$ along 
$H^\perp$ is of odd order. But by hypothesis, this eigenvalue is a root of unity belonging to the field of definition of $W$.
Its order is then a power of $2$. Hence $w$ acts trivially on $H^\perp$ which means that $w \in C_H$. 

If $G$ is a $2$-subgroup of $W$ then every element in $p^{-1}(G)$ has infinite order thanks to Corollary~\ref{cor-lifting}.
Let $P_0=p^{-1}(Z(W)) \isom{gr.} \ZZ^{|\Hscr|}$ (see~\cite{marin2}), the extension $\extensionzu{P_0}{p^{-1}(G)}{\barre{G}}{}{}$ gives the result
(let us remind that $\barre{G}$ is the direct image of $G$ in $W/ZW$).  

With Theorem~\ref{theo-finite-order}, we get the converse. If $G$ is not a $2$-group, then it contains an element of odd order.
The first part ensures us that this element may be lifted in $p^{-1}(G)$ as a finite order element.

Let us now consider $G \subset W$ whose order is odd. Then, for every $H \in \Hscr$ and every $w \in G\cap N_H$, we have
$w \in C_H$ by the above argument since the order of $w$ is odd. Hence~Theorem~\ref{theo-finite-order} gives the result.   
\end{demo}

\begin{rem} The irreducible complex reflection groups $W$ such that every element of odd order has a finite order lifting 
in $B/[P,P]$ are the same as the list of Corollary~\ref{cor-real-reflection-group} that is to say
$G_8, G_9, G_{12}, G_{13}, G_{22}, G_{23}, G_{24}, G_{28}$, $G_{29}, G_{30}, G_{31}, G_{35}, G_{36}, G_{37}$ 
for the exceptional ones, and $G(d,1,1)$ with $d$ a power of $2$, 
$G(de,e,r)$ with $d$ and $e$ powers of $2$ and $r \geq 2$, $\Sfr_r$ and $G(e,e,2)$ with $e \geq 3$. 

The preceding corollary shows that all these groups verify the property. 
For the exceptional ones, a~\cite{gap} computation using \cite{chevie} gives that there are no more exceptional group such that every element of odd order has a finite order lifting in $B/[P,P]$. 
Let us now study the infinite series and consider $G(de,e,r)$ such that every element of odd order has a finite
order lifting in $B/[P,P]$. If $r=1$, since $B/[P,P]=\ZZ$ then $G(d,1,1)$ should have only one odd order element, so that
$d$ is a power of $2$. Let us assume that $r \geq 2$ and $d \geq 2$. For every $\zeta \in \UU_{de} \setminus \{1\}$, the element 
$\diag(\zeta,\zeta^{-1},1,\ldots, 1)$ does not have a finite order lifting in $B/[P,P]$ (see Corollary~\ref{cor-infinite-series}). 
It implies that $de$ has to be a power of $2$. Let us consider the case where $r \geq 3$ and $d=1$.
For $\zeta \in \UU_{de} \setminus \{1\}$, the element $\diag(\zeta,\zeta,\zeta^{-2},1,\ldots,1)$ does not have a finite order lifting
in $B/[P,P]$ (see Corollary~\ref{cor-infinite-series}). It implies that $de$ has to be a power of $2$. 
When $r=2$ et $d=1$, we get $G(e,e,2)$ which verify the property. 
\end{rem}

\begin{rem} Following Cayley's argument, every odd-order group may be identified to a subgroup of $\Sfr_n$ for $n$ large enough
and thus to a subgroup of $B_n/[P_n,P_n]$ thanks to Corollary~\ref{cor-real-reflection-group}. In the course of writing, 
Gonçalves, Guaschi and Ocampo have communicated to us that they independently got this same result. 
\end{rem}

\subsection{Infinite Braid Group}

Let $B_\infty$ the infinite braid group. This is the direct limit of the family $(B_n)_{n \in \NN}$ 
where $B_n$ is the standard braid group on $n$ strands and the map $B_n \rightarrow B_{n+1}$ consists to add one strand on the right. 
The pure braid group $P_\infty$ is the direct limit
of the family $(P_n)_{n \in \NN}$ where $P_n$ is the pure braid group on $n$ strands. 

Since $[P_{n-1},P_{n-1}]=[P_n,P_n] \cap B_n$,  $B_{n-1}/[P_{n-1},P_{n-1}]$ may be identified to a subgroup 
of $B_n/[P_n,P_n]$ and then also to a subgroup of the direct limit of the family $(B_n/[P_n,P_n])_{n \in \NN}$ which 
is nothing else than $B_{\infty}/[P_\infty,P_\infty]$. Moreover, with the same kind of arguments 
$P_{\infty}/[P_\infty,P_\infty]$ is the direct limit of the groups $P_n/[P_n,P_n]$ and is a free abelian group 
on the set of $2$-sets of the infinite set $\NN^*$. Finally, the group $\Sfr_\infty$ may be defined as 
the direct limit of the family $(\Sfr_n)_{n\in\NN}$ where an element of $\Sfr_{n-1}$ is seen as an element of $\Sfr_n$ fixing $n$.
But $\Sfr_\infty$ may also be viewed as the group of permutation of $\NN^*$ whose support is finite.

\begin{prop} For every odd-order group $G$, there exists a subgroup of 
$B_{\infty}/[P_\infty,P_\infty]$ isomorphic to $G$ and $B_{\infty}/[P_\infty,P_\infty]$ contains no even order element. 
\end{prop}

\begin{demo} The group $B_{\infty}/[P_\infty,P_\infty]$ is the direct limit of 
the family $(B_n/[P_n,P_n])_{n \in \NN}$. 
Moreover Proposition~2.4 of~\cite{marin2} ensures us that the maps
$B_{n-1}/[P_{n-1},P_{n-1}] \rightarrow B_n/[P_n,P_n]$ are injective for all positive $n$. 
Then the map $B_n/[P_n,P_n] \rightarrow B_\infty/[P_\infty,P_\infty]$ is injective for all $n$. 

Let us now consider $G$ a group of odd order. Following Cayley's argument, we let $G$ acts on itself
by left translation, so that $G$ identifies to an odd-order group of $\Sfr_{|G|}$ the symmetric group on $|G|$ letters. 
But Corollary~\ref{cor-real-reflection-group} ensures us that $G$ identifies 
to a subgroup of $B_{|G|}/[P_{|G|},P_{|G|}]$ which itself is a subgroup of $B_\infty/[P_\infty,P_\infty]$.

Let $x$ be an even order element of $B_{\infty}/[P_\infty,P_\infty]$, then it belongs to $B_n/[P_n,P_n]$ for some $n$. 
But $B_n/[P_n,P_n]$ has no even order element. 
\end{demo}

We can be even more precise. 

\begin{prop}\label{prop-subgroup-binfty} 
Let $\Sfr_\infty$ be the permutation group of the infinite set $\NN^*=\{1,2,\ldots\}$ with finite support. 
We have the following extension
$\extensionuu{{P_\infty}^\textrm{ab}}{B_\infty/[P_\infty,P_\infty]}{\Sfr_\infty}{}{p}$ and
${P_\infty}^\textrm{ab}$ identifies to the free abelian group over the set $\Pscr_2$ of $2$-subsets of $\NN^*$. 

A finite subgroup of $B_\infty/[P_\infty,P_\infty]$ is of odd order and maps isomorphically to $\Sfr_\infty$ through $p$. 
For every odd-order subgroup $G$ of $\Sfr_\infty$, there exists a group homomorphism 
$s:G \ra B_\infty/[P_\infty,P_\infty]$ such that $ps=\id_G$.
Moreover two such group homomorphisms are conjugated by an element of ${P_\infty}^\textrm{ab}$ and the intersection of 
the normalizer of a finite subgroup $\wt{G}$ of $B_{\infty}/[P_\infty,P_\infty]$ with ${P_\infty}^\textrm{ab}$ 
is isomorphic to $\ZZ \Pscr_2/p(\wt{G})$.
\end{prop}

\begin{demo} A finite subgroup of $B_\infty/[P_\infty,P_\infty]$ is contained in some 
$B_n/[P_n,P_n]$ and hence of odd order~(Corollary~\ref{cor-real-reflection-group}). Moreover 
since $P_\infty/[P_\infty,P_\infty]$ is a free abelian group, the surjective map $p$ preserves the order 
of the finite order elements (we have already seen this in the proof of theorem~\ref{theo-finite-order}). 
So a finite subgroup of $B_\infty/[P_\infty,P_\infty]$ maps isomorphically onto a subgroup of $\Sfr_\infty$ through
$p$. 

Let $G$ be a finite subgroup of $\Sfr_\infty$. Then it is a subgroup of $\Sfr_n$ for some $n$ and hence embeds in $B_n/[P_n,P_n]$
and so in $B_\infty/[P_\infty,P_\infty]$. The second part follows from the fact that $H^1(G,P_\infty/[P_\infty,P_\infty])=\{0\}$
since $P_\infty/[P_\infty,P_\infty]$ is a permutation module (Lemma~\ref{lem-ext}). 
Finally, let $x \in {P_\infty}^\textrm{ab}$. Then, 
as in the proof of Theorem~\ref{theo-finite-order}, $x$ verifies $x\wt{G}x^{-1}=\wt{G}$ if and only if $gxg^{-1}=x$ 
for every $g \in \wt{G}$.
\end{demo}	
	
\subsection{Some More Examples}

We have seen that every odd-order subgroup of $\Sfr_n$ can be lifted in $B_n/[P_n,P_n]$.
This has been applied to show that every odd-order group $G$ can be embedded in $B_\infty/[P_\infty,P_\infty]$. 
Here we describe some odd-order subgroups of finite symmetric groups verifying a stronger property: they do not 
meet any stabilizer of a hyperplane. 

The study of such subgroups is quite natural because they are the subgroups of $W$ that acts freely on $\Hscr$.
Moreover, when $G$ is a subgroup of $W$ such that for every $H \in \Hscr$, 
$G \cap N_H = \{1\}$ then the criterion of~Theorem \ref{theo-finite-order} is clearly verified for $G$. 
But we have in fact a stronger property:
$H^2(G,\ZZ\Hscr)$ is trivial. In order to prove this, it suffices to consider the orbits of $\Hscr$ under $G$ and to apply
Shapiro's isomorphism.

\subsubsection{The case of the symmetric group}

To study the case of the symmetric group, let us introduce the following notation. 
For an integer $n$, we denote by $F_n$ the set 
$$F_n=\{\si \in \Sfr_n, \quad \si \textrm{ is of cycle type } k^{[n/k]} \textrm{ or } 1^{[1]}k^{[(n-1)/k]} 
\textrm{ with } k \textrm{ odd}\}\,.$$

\begin{prop} Let $G$ be a subgroup of $\Sfr_n$. For $i \neq j$, we denote by $C(i,j)$ the centralizer of the transposition 
$(i,j)$: this is also the stabilizer $N_{H_{i,j}}$. Then, $G \cap C(i,j)= \{\id\}$ for all $i\neq j$ if and only if $G \subset F_n$.
\end{prop}

\begin{demo} We have $C_{i,j}=\Sfr_{n-2}\times \langle (i,j) \rangle$. In particular, $w \in C_{i,j}$ if and only if 
$\{w(i),w(j)\}=\{i,j\}$. 

Let now $G$ such that $G\cap C(i,j)=\{\id\}$ for all $i \neq j$. If $G$ has an element $w$ of order $n=2k$ then 
$w^k$ is a non trivial product of transpositions and then belongs to $C(i,j)$ for some $(i,j)$. So 
every element of $G$ is of odd order and every cycle of an element of $G$ is of odd length 
(to see this, we could also have applied the fact the $G$ lifts into $B/[P,P]$ thanks to Theorem~\ref{theo-finite-order}). 
Moreover every non trivial element $w$ of $G$ has at most one fixed point: if there are 
two or more, then let $\{i,j\}$ be two such fixed points, then $w \in C(i,j)$. 
We then deduce that the length of every non trivial cycle of $w \in G$ is the same. Indeed if $c_1$ and $c_2$
are two non trivial cycles of length $k_1 < k_2$ then $w^{k_1} \neq 1$ and has more than one fixed point. 
This gives the wanted cycle decomposition. 

Conversely, assume that $G \subset F_n$. For $g \in F_n$, every power of $g$ is still in $F_n$ since 
every cycle of $g$ will decomposes in cycles of length a divisor of an odd number. 
So for every $i\neq j$, the set $\{i,j\}$ is not stable by any non trivial element of $G$.
\end{demo}

%
%
%

\begin{ex}[Group of odd order] Let $G$ be a group of odd order. When $G$ acts on itself by left multiplication, 
then the cycle decomposition of an element $g \in G \subset \Sfr_{|G|}$ is a product of $[G:\langle g \rangle]$
cycles of length the order of $g$ which is odd. And so $G \subset F_{|G|}$ which a bit more precise that what was needed 
in~Proposition~\ref{prop-subgroup-binfty}.
\end{ex}

Let us know study a bit Frobenius groups who will provide some more examples. 
Let $G$ be a Frobenius group with kernel $N$ and $H$ a Frobenius complement. The group $G$ is a disjoint union
$$G= N \sqcup \bigsqcup_{g \in G/H} (gHg^{-1}\setminus {1})\,.$$

\begin{lem}\label{lem-action-frobenius-group}
The group $G$ acts on $G/H$ with the following properties. 

\begin{ri}{(ii)}
\item An element $n \in N \setminus \{1\}$ has no fixed point on $G/H$. 
\item For $g \in G$, an element of $gHg^{-1}\setminus {1}$ has only one fixed point on $G/H$ namely $gH$.  
\end{ri}
\ni In particular, $G$ acts faithfully on $G/H$. 
\end{lem}

\begin{demo} 
\begin{ri}{(ii)} 
\item If $ngH=gH$ then $g^{-1}ng \in H$ and $n \in gHg^{-1}$ which is absurd thanks to the set partition of $G$. 
\item If $ghg^{-1}g'H=g'H$ (with $h \in H \setminus \{1\}$) then $g'^{-1}ghg^{-1}g' \in H$ and 
$h \in H \cap g^{-1}g'Hg'^{-1}g$. Hence, the set partition of $G$ ensures us that $g^{-1}g' \in H$.
\end{ri}
\end{demo}

\begin{cor} For every $g \in G$, the cycle decomposition of $g \in G$ seen as a permutation of $G/H$ is the following
\begin{ri}{(ii)}
\item if $g \in N$, then $g$ is a product of $[G:H]/k$ cycles of length $k$ where $k$ is the order of $g$.  

\item if $g \in g'Hg'^{-1} \setminus \{1\}$ then $g$ is a product of $([G:H]-1)/k$ cycles of length $k$ where $k$ is the order of $g$. 
\end{ri}
\end{cor}

\begin{demo} It is trivial when $g =1$. Let us consider $g \neq 1$.
\begin{ri}{(ii)}
\item For $g'H \in G/H$ and $1 \leq \ell <k$, we have $g^\ell g'H \neq g'H$ since $g^\ell \in N \setminus \{1\}$.
Hence the orbit of $g'H$ under $\langle g \rangle$ has $k$ elements. 

\item Let $1 \leq \ell < k$ and $g''H \in G/H$ with $g''H \neq g'H$. 
We have $g^\ell g''H\neq g''H$ since $g^\ell \in g'Hg'^{-1} \setminus \{1\}$
and $g''H \neq g'H$. Hence the orbit of $g''H$ under $\langle g \rangle$ has $k$ elements. 
\end{ri}
\end{demo}

\begin{cor} Let $K$ be a subgroup of a Frobenius group of odd order, then $K$ is contained in $F_{(G:H)}$.
\end{cor}

This last corollary associated to Corollary~\ref{cor-real-reflection-group} generalizes Corollary~3.10 and~3.11 
of~\cite{marin2} and Theorem~7 of~\cite{ggo}.

\subsubsection{The infinite series}

Let us introduce the following subset of $G(de,e,r)$
$$F(de,e,r)=\{w \in G(de,e,r), \textrm{the cycle type of } \si_w  \textrm{ is } [k]^{r/k} \textrm{ with } k \textrm{ odd and } 
	p_{w,\si}=1 \textrm{ for every } \si \textrm{ cycle of } \si_w\}$$
\ni where $\si_w$ is permutation associated to $w$. 

\begin{prop} Let us consider $G$ a subgroup of $G(de,e,r)$ with $d \geq 2$ and $\Hscr$ be the set of hyperplanes of $G(de,e,r)$.  
Then $G \cap N_H= \{1\}$ for every $H \in \Hscr$ if and only if $G \subset F(de,e,r)$. 
\end{prop}

\begin{demo} Assume that $G \cap N_H= \{1\}$ for every $H \in \Hscr$ then from Theorem~\ref{theo-finite-order} we 
deduce that every element of $G$ is of odd order and so for every $w \in G$ and every $\si$ cycle of $\si_w$ 
is of odd order and from Corollary~\ref{cor-infinite-series} that $p_{w,\si}=1$
for every $\si$ cycle of $\si_w$ and every $w \in G$. Moreover consider $w \in G$ such that $\si_w$ has at least 
two cycles of different lengths. Suppose that these lengths are $\ell_1 < \ell_2$. 
Then $w^{\ell_1} \in G \setminus \{1\}$ and stabilizes every hyperplane of the form $\{z_i=0\}$ where $i$ belongs to the  
orbit of length $\ell_1$. 

Conversely, assume that $G \subset F(de,e,r)$. If $w \in G$ stabilizes the hyperplane $\{z_i=0\}$
then $\si_w(i)=i$ and hence $\si_w=\id$ since every cycle of $\si_w$ has the same length. 
Since $p_{w,\si}=1$ for every $\si$ cycle of $\si_w=1$, we get that $w=1$. 
If $w \in G$ stabilizes $H_{i,j,\zeta}$ for $i \neq j$ and $\zeta \in \UU_{de}$ then $\{\si_w(i),\si_w(j)\}=\{i,j\}$. 
Since every cycle of $\si_w$ is of odd length then $\si_w(i)=i$ and $\si_w(j)=j$. Hence $w$ stabilizes $\{z_i=0\}$
so $w=1$. 
\end{demo}

%


\bibliographystyle{mythese2}
    \begin{small}
        \bibliography{mythese-admin-utf8}
    \end{small}

\end{document}